\documentclass[10pt]{article}

\usepackage{graphicx,hyperref}
\usepackage{amsmath}
\usepackage[version=4]{mhchem}
\usepackage{siunitx}
\usepackage{longtable,tabularx}
\usepackage{pstricks}
\usepackage{amsfonts}
\usepackage{latexsym}
\usepackage{enumerate}
\usepackage{epsfig}
\usepackage{subfig}
\usepackage{framed,fancybox}
\usepackage{bbm}
\usepackage{multirow}
\usepackage[T1]{fontenc}
\usepackage{hhline}

\usepackage{threeparttable}
\usepackage{rotating}
\usepackage{amssymb,amscd}
\makeatletter
\def\myVCENTER#1{\vcenter{\hbox{$\m@th#1$}}}
\makeatother

\long\def\symbolfootnote[#1]#2{\begingroup\def\thefootnote{\fnsymbol{footnote}}\footnote[#1]{#2}\endgroup}

\definecolor{shadecolor}{gray}{0.99}
{\endMakeFramed}

\newenvironment{shadedframe}{%
 \MakeFramed {\FrameRestore}}
{\endMakeFramed}

\definecolor{shadecolor}{gray}{0.99}
\long\def\symbolfootnote[#1]#2{\begingroup\def\thefootnote{\fnsymbol{footnote}}\footnote[#1]{#2}\endgroup}
\def\qed{\hfill{$\vcenter{\hrule height1pt \hbox{\vrule width1pt height5pt
    \kern5pt \vrule width1pt} \hrule height1pt}$} \medskip}
\newcommand{\m}[1]{{\bf{#1}}}
\newcommand{\g}[1]{\boldsymbol #1}
\newcommand{\bb}[1]{\mathbb #1}

\newcommand{\C}[1]{{\cal {#1}}}

\DeclareMathOperator{\erf}{erf}

\title{\bf Method for Solving Chance Constrained Optimal Control \\ Problems Using Biased Kernel Density Estimators}

\author{Rachel E.~Keil\thanks{Ph.D.~Candidate, Department of Mechanical and Aerospace Engineering, University of Florida, Gainesville, FL 32611-6250. E-mail: rekeil@ufl.edu} \\ Alexander T.~Miller \thanks{Ph.D.~Candidate, Department of Mechanical and Aerospace Engineering, University of Florida, Gainesville, FL 32611-6250. E-mail: alexandertmiller@ufl.edu}
 \\ Mrinal Kumar \thanks{Associate Professor, Department of Mechanical and Aerospace Engineering, The Ohio State University, Columbus, OH 43210. AIAA Senior Member. E-mail: kumar.672@osu.edu} \\ Anil V.~Rao\thanks{Professor, Erich Farber Faculty Fellow and University Term Professor, Department of Mechanical and Aerospace Engineering, University of Florida, Gainesville, FL 32611-6250.  E-mail: anilvrao@ufl.edu.  Corresponding Author.}}
\date{}
\begin{document}
\maketitle
\thispagestyle{empty}

\begin{abstract}
A method is developed to numerically solve chance constrained optimal control problems.  The chance constraints are reformulated as nonlinear constraints that retain the probability properties of the original constraint.  The reformulation transforms the chance constrained optimal control problem into a deterministic optimal control problem that can be solved numerically.  The new method developed in this paper approximates the chance constraints using Markov Chain Monte Carlo (MCMC) sampling and kernel density estimators whose kernels have integral functions that bound the indicator function.  The nonlinear constraints resulting from the application of kernel density estimators are designed with bounds that do not violate the bounds of the original chance constraint.  The method is tested on a non-trivial chance constrained modification of a soft lunar landing optimal control problem and the results are compared with results obtained using a conservative deterministic formulation of the optimal control problem.  The results show that this new method efficiently solves chance constrained optimal control problems.
\end{abstract}


\renewcommand{\baselinestretch}{2}
\normalsize\normalfont 

\section{Introduction}
Optimal control problems arise frequently in different fields including various branches of engineering and disciplines outside of engineering.  The goal of solving an optimal control problem is to determine the state and control of a controlled dynamic system that optimizes a performance index subject to dynamic constraints, path constraints and boundary conditions~\cite{Betts3}.  The applications for optimal control problems include robotics, disease control, chemical processes, and flight.   Optimal control problems can be broadly categorized as either deterministic or stochastic.  In a deterministic optimal control problem, no part of the optimal control problem contains uncertainty.  By contrast, in a stochastic optimal control problem uncertainty can be present in any part of the optimal control problem.  Uncertainty in an optimal control problem arises in various forms including process noise, measurement noise, and uncertainty in the constraints.  When uncertainty is present in the constraints of the stochastic optimal control problem, a particular class of these optimal control problems are known as chance constrained optimal control problems (CCOCPs).   

Obtaining a solution to a CCOCP analytically is possible only for simple problems.  Analytical treatment of chance constraints has been most extensively explored in the domain of robust model predictive control (RMPC), with applications that include traffic control~\cite{yan05}, chemical process control~\cite{li08}, and others~\cite{li02, farina13}.  Due to the small number of problems solvable analytically, solutions to most CCOCPs must be obtained using numerical methods.  Numerically solving a CCOCP is challenging due to the probabilistic formulation of the chance constraints.  A simple approach to solving a CCOCP numerically is to reformulate the chance constraints as worst case scenario deterministic constraints.  This reformulation transforms the CCOCP to a deterministic optimal control problem.  It is noted, however, that this worst case scenario deterministic optimal control problem results in a very conservative solution.  Additionally, this deterministic optimal control problem could even prove to be infeasible~\cite{Keil1}.  The issues with the worst case scenario method are in large part due to not incorporating stochastic properties in the deterministic constraints.  As a result, it would be preferable to reformulate the chance constraints as deterministic constraints in an analogous way to how stochastic processes are often reformulated as deterministic processes that incorporate key stochastic properties (for example, the key properties of a Gaussian distribution as used in a Kalman filter).  If such a formulation could be obtained for the chance constraints, the CCOCP would be transformed to a numerically solvable deterministic optimal control problem that incorporates key properties of the CCOCP.  

Various methods have been developed to transform chance constraints to nonlinear constraints that retain the main stochastic characteristics of the chance constraints~\cite{geletu13}.  In the method of Ref.~\cite{blackmore10}, weighted summations are used to approximate the chance constraints as deterministic constraints.  Next, in the methods of Refs.~\cite{Blackmore1,ono10,okamoto19}, the chance constraint is transformed using Gaussian properties to a deterministic form.  Note, however, that the methods of Refs.~\cite{blackmore10,Blackmore1,ono10,okamoto19} are applicable only to linear chance constrained optimization problems.  In addition, Ref.~\cite{hokayem13} transforms the chance constraint to a deterministic form using Gaussian properties.  The approach of Ref.~\cite{hokayem13} can be applied to chance constrained linear quadratic Gaussian problems.   In contrast, the method of Ref.~\cite{muhlpfordt18} uses arbitrary non-Gaussian distributions to determine a deterministic approximation of the chance constraint.  It is noted, however, that the method of Ref.~\cite{muhlpfordt18} is only applicable to linear chance constrained optimization problems.  A method applicable to a more general class of chance constrained optimization problems is presented in Ref.~\cite{Pinter}, which provides various deterministic expressions to approximate a chance constraint based on available information about the chance constraint.  The method of Ref.~\cite{Pinter} can be difficult to implement depending on the form of the constraint.  In Ref.~\cite{Nemirovski}, under the assumptions of convex and affine relations for the chance constraint, a convex chance constraint approximation is constructed.  As such, the method of Ref.~\cite{Nemirovski} is also applicable only to a specific subset of chance constraints.   Another method of interest is the sample average approximation method of Ref.~\cite{Pagnoncelli1}.  In the method of Ref.~\cite{Pagnoncelli1}, an approximate deterministic constraint is applied using the indicator function, a function form which can increase the numerical difficulty of obtaining a solution.  In the method of Ref.~\cite{ono15}, linear approximations that depend on the indicator function are employed to effectively transform the chance constrained optimization problem to a deterministic optimization problem.  This dependence on the indicator function can increase computational expense in solving this determinstic optimization problem.  A method that is more computationally efficient than that of Ref.~\cite{ono15} is the scenario approach of Refs.~\cite{Calafiore1,Calafiore2}.  The scenario approach determines the worst case scenario objective by using deterministic constraints resulting from the application of a finite number of samples of the uncertainty to the chance constraints.  The disadvantage of the scenario approach is that without numerically restrictive fine tuning~\cite{Campi1}, the resulting solutions are overly-conservative.  This conservatism results from the nonlinear constraint approximation being conservative relative to the original chance constraint.  As such, more recent methods to transform chance constraints to nonlinear constraints are designed to be less conservative as well as being applicable to a wider range of chance constraints.

More recently, sampling methods have been developed for transforming chance constraints to not overly-conservative nonlinear constraints.  Similar to the scenario approach, these methods employ sampling to approximate the chance constraint.  In the method of Refs.~\cite{Chai,Kumar1,Ahmed}, samples are applied to obtain an approximate expectation of various different functions that is a nonlinear approximation of the chance constraint which provides an upper bound on the chance constraint.  The disadvantage of the method of Refs.~\cite{Chai,Kumar1,Ahmed} is that the resulting nonlinear approximation can still be too conservative relative to the chance constraint.  Another class of methods that applies sampling and kernels leads to less conservative nonlinear constraints.  In one such method from Ref.~\cite{Bharath}, kernels are used with the samples in Reproducing Kernel Hilbert Space to approximate the chance constraint.  This method can become computationally difficult as the dimensions of the problem grow.  Another method applies kernel density estimators (KDEs) to the samples to obtain a nonlinear approximation of the chance constraint~\cite{Calfa,Caillau}.  KDEs result in a nonlinear constraint that is not overly-conservative relative to the chance constraint, at the cost of potentially violating the bounds of the chance constraint unlike the method of Refs.~\cite{Chai,Kumar1,Ahmed}.  Developing a method that combines the advantages of the method from Refs.~\cite{Chai,Kumar1,Ahmed} and the KDE method from Refs.~\cite{Calfa,Caillau} would be an improvement on previously developed methods.

In this paper, the criteria for a kernel to provide a KDE that does not violate the bounds of the chance constraint is developed using key results from Refs.~\cite{Chai,Kumar1,Ahmed}.  This criteria will be used to define a new method to approximate the chance constraints as KDEs that provide an upper bound on the chance constraint, while not being overly-conservative relative to the chance constraint.   In particular, the new method transforms chance constraints to deterministic nonlinear constraints.  Consequently, the CCOCP is transformed to a deterministic optimal control problem.  Numerical methods can then be applied to solve the transformed CCOCP.  

Various numerical methods have been developed for solving deterministic optimal control problems.  These methods can be categorized broadly as either indirect or direct methods~\cite{Stryk0}.  Indirect methods use the first order optimality conditions to transform the optimal control problem to a Hamiltonian boundary value problem (HBVP) that is then solved.  For direct methods the original optimal control problem is directly solved by parameterizing the states and possibly the control.  The direct method of Gaussian quadrature collocation in the form of either Legendre-Gauss (LG)~\cite{Benson2,Rao8}, Legendre-Gauss-Lobatto (LGL)~\cite{Elnagar1} or Legendre-Gauss-Radau (LGR)~\cite{Garg1,Garg2,Patterson2015} collocation has broad applicability, in addition to providing high accuracy solutions with exponential convergence~\cite{Canuto1,Fornberg1}.  Direct collocation can also be easily paired with readily available software~\cite{Gill1,Biegler1}.  As a result of the advantages of direct Gaussian quadrature collocation, in this paper LGR collocation will be applied with the new method for transforming chance constraints to deterministic nonlinear constraints to develop a computational framework for efficiently solving CCOCPs in their transformed deterministic form.    

The paper is organized as follows.  Section~\ref{sect:CCnonlinear} gives a brief overview of methods applied to transform chance constraints to nonlinear constraints and uses aspects of these methods to develop the new method.  Section~\ref{sect:bolza} describes a general chance constrained optimal control problem.  Section~\ref{sect:LGR} provides a brief overview of the LGR collocation method.  Section~\ref{sect:method} describes the computational framework to numerically solve the transformed CCOCP using available optimal control software.  In Section~\ref{sect:examples}, this computational framework is applied to an example.  Section~\ref{sect:conc} provides some conclusions.

\section{Chance Constrained Optimization}\label{sect:CCnonlinear}
In this section, a new method is developed for reformulating a chance constrained optimization problem as a deterministic optimization problem.  The reformulation is accomplished by transforming the chance constraints to nonlinear constraints that retain key stochastic properties of the chance constraints.  The new method developed combines two previously developed methods, and these two methods are described in Section~\ref{sect:previous}.  The first of the methods, Method 1 described in Section~\ref{sect:firstm}, defines a nonlinear approximation of the chance constraint.  This approximation is an approximation of the cumulative distribution function (CDF) of the random variable associated with the chance constraint.  The advantage of Method 1 is that the approximation is an upper bound on the chance constraint.  The disadvantage of Method 1 is that this approximation of the CDF is allowed to exceed unity.  The second of the methods, Method 2 described in Section~\ref{sect:secondm}, uses Kernel Density Estimators (KDE)s to obtain a nonlinear approximation of the chance constraint.  The advantages of Method 2 are that the approximation of the CDF stemming from the KDE does not exceed unity and that there is is a wide range of functions that can be used for a KDE.  The disadvantage of Method 2 is that the approximation may violate the bound on the chance constraint.  Consequently, a new method that combines the advantages of these previously developed methods would be an improvement on said methods.  

One particular nonlinear approximation used for Method 1 is the Split-Bernstein approximation described in Ref.~\cite{Kumar1}.  This nonlinear approximation is shown in Section~\ref{sect:SBkernel} to be a KDE with the appropriate choice of parameters.  The Split-Bernstein approximation with certain parameters combines the advantage of Method 1 in upper bounding the chance constraint with one of the advantages of Method 2 in having an approximation of the CDF that does not exceed unity.  As a result of this last fact, in Section~\ref{sect:BiasKDE} the Split-Bernstein KDE is analyzed to determine the necessary criteria for other KDEs to retain these advantages.  In the new method, a nonlinear approximation to the chance constraint is defined that applies this necessary criteria in order to use KDEs that ensure an upper bound the chance constraint.  This new method is then applied to chance constrained optimal control, as described in Section~\ref{sect:bolza}.  

\subsection{Previous Methods for Reformulating Chance Constraints}\label{sect:previous}
To begin the development of the new method, consider the following general chance constrained optimization problem: 
\begin{equation}\label{eq:CCexample_no}
\begin{array}{c}
 \min\limits_{z \in Z} J(z) \\
 P(\m{F}(\m{z},\g{\xi}) \geq \m{q}) \geq  1-\epsilon,
\end{array}
\end{equation}
where $\m{z}$ is a decision variable defined on the feasible set $\m{Z} \subset \bb{R}^n$ and $\g{\xi}$ is a random vector supported on set $\Omega \subseteq \bb{R}^d$ for which samples of $\g{\xi}$ will be used.  The function $ \m{F}(\m{z},\xi) \geq \m{q}$ is an event in the probability space $P(\cdot)$, where $\m{F}( \cdot )$ maps $\bb{R}^n \times \bb{R}^d \rightarrow \bb{R}^{n_g}$ and $\epsilon$ is the risk violation parameter.  
Because the function $\m{F}(\m{z},\g{\xi}) $ is a vector, Eq.~\eqref{eq:CCexample_no} is a joint chance constraint.  Using Boole's inequality together with the approach of Refs.~\cite{Blackmore1} and~\cite{Kumar1}, the chance constraint given in Eq.~\eqref{eq:CCexample_no} can be redefined in terms of the following two conservative constraints (see Refs.~\cite{Blackmore1} and~\cite{Kumar1} for the proof):
\begin{equation}\label{eq:CC_scalar}
\begin{array}{c}
P( \psi \geq q_m  ) \geq  1-\epsilon_m, \\
\sum\limits_{m =1}^{n_g} \epsilon_m \leq \epsilon,
\end{array}
\end{equation} 
or, equivalently 
\begin{equation}\label{eq:eventcomp}
\begin{array}{c}
P( \psi < q_m  ) \leq  \epsilon_m, \\
\sum\limits_{m =1}^{n_g} \epsilon_m \leq \epsilon,
\end{array}
\end{equation}
where $m\in [1,\dots,n_g]$ is the index corresponding to the $m$th component of the event and $\psi$ is a random variable with a PDF that is assumed to be unknown, defined such that
\begin{equation}
\psi = F_m(\m{z},\g{\xi}).
\end{equation}  
It should be noted that in Ref.~\cite{Nemirovski} the first of the two scalar constraints in Eq.~\eqref{eq:CC_scalar} was also applied, but not the second as this constraint would destroy the convexity of the problem.  As the retention of convexity is not required for application of numerical methods, both constraints can be used to replace the chance constraint of Eq.~\eqref{eq:CCexample_no}.

\subsubsection{Method 1: Nonlinear Approximation of the Chance Constraint}\label{sect:firstm}
In the first method, an upper bound on the chance constraint is defined as the expectation of some function $w(\g{\xi})$.  Because the expectation cannot be obtained analytically, a nonlinear approximation of the expectation is obtained using sampling.  In particular, samples of $\g{\xi}$ are obtained using Markov Chain Monte Carlo (MCMC) sampling or by some simpler method when advantageous properties of $\g{\xi}$ are assumed.  Formally, the upper bound on the chance constraint given in Eq.~\eqref{eq:CC_scalar} for the first method is defined in terms of the expectation $\bb{E}$ and the indicator function $1_{( \cdot )}$ as:
\begin{equation}\label{eq:equalApp}
1-\epsilon_m \ \leq P(\psi \geq q_m) = \bb{E}_{\g{\xi}} [	1_{[q_m,+\infty)}(\psi-q_m)] \leq \bb{E}_{\g{\xi}} [w( \psi-q_m)],
\end{equation}
where $w(\cdot)$ represents a function that has the property
\begin{equation}\label{eq:boundind}
1_{[q_m,+\infty)}(\psi-q_m) \leq w( \psi-q_m), \ \forall \left( \psi-q_m \right).
\end{equation}
Different nonlinear approximations of the upper bound on the chance constraint from Eq.~\eqref{eq:equalApp} were developed in Refs.~\cite{Chai,Kumar1,Ahmed}.  The property of the function $w(\cdot)$ in Eq.~\eqref{eq:boundind} and the type of samples obtained for $\g{\xi}$ are what ensure that the nonlinear approximation of the expectation still provides an upper bound on the chance constraint.  The function $w(\cdot)$, however, has the disadvantage of being allowed to pass unity instead of being bounded by it.  As a result, the nonlinear approximation of the chance constraint of Eq.~\eqref{eq:equalApp} resulting from Method 1 will exceed unity and so cannot be an approximate CDF.  Consequently, in the next section a second method that uses KDEs to obtain a nonlinear approximation of the chance constraint which is an approximate CDF is considered.

\subsubsection{Method 2: Kernel Density Estimator}\label{sect:secondm}
The second method, described in Ref.~\cite{Caillau}, uses KDEs to obtain a nonlinear approximation of the chance constraint.  The nonlinear approximation is obtained by first applying kernels to determine an approximate PDF $\hat{f}_{\psi}$ of a desired random variable $\psi$, using samples $\psi_j$ of the random variable.  The approximate PDF is defined as follows:
\begin{equation}\label{eq:KDE}
\hat{f}_{\psi}(q_m) = \frac{1}{N} \sum_{j = 1}^N \frac{1}{h} k \left( \eta_j \right),
\end{equation}
where $k(\cdot)$ is the kernel, $h$ is the bandwidth, $j = 1,\dots,N$ is the number of samples and $\eta_j$ is defined as
\begin{equation}\label{eq:etadef}
 \eta_j = \frac{q_m-\psi_j}{h}.
\end{equation}
The right hand side of Eq.~\eqref{eq:KDE} is the KDE.  To define the relation between the chance constraint and the KDE, first the approximate probability $\hat{p}$ is defined as
\begin{equation}\label{eq:CDFproof}
\hat{p}(\psi < q_m) = \int_{-\infty}^{q_m} \hat{f}_{\psi}(x)d x.
\end{equation}
Then, applying the general form of a KDE from Eq.~\eqref{eq:KDE} to the right hand side of Eq.~\eqref{eq:CDFproof}:
\begin{equation}\label{eq:Probsum}
\int_{-\infty}^{q_m} \hat{f}_{\psi}(x)d x = \frac{1}{Nh}\sum_{j = 1}^N \int_{-\infty}^{q_m} k \left( \frac{x-\psi_j}{h} \right) d x.
\end{equation}
Using the definition of $\eta_j$ from Eq.~\eqref{eq:etadef}, Eq.~\eqref{eq:Probsum} can be reformulated as:
\begin{equation}\label{eq:CDFactual}
\frac{1}{Nh}\sum_{j = 1}^N \int_{-\infty}^{q_m} k \left( \frac{x-\psi_j}{h} \right) d x = \frac{1}{N} \sum_{j = 1}^N \int_{-\infty}^{\eta_j}k(v_j) d v_j.
\end{equation}
The following integrated function of the kernel can now be defined as:
\begin{equation}\label{eq:kernCDF}
K(\eta_j) = \int_{-\infty}^{\eta_j}k(v_j) d v_j.
\end{equation}
Using Eq.~\eqref{eq:kernCDF}, the approximation of the chance constraint given in Eq.~\eqref{eq:CDFproof} can then be defined relative to a general KDE as 
\begin{equation}\label{eq:CDF}
\hat{p}(\psi < q_m) = \frac{1}{N} \sum_{j = 1}^N K(\eta_j),
\end{equation}
where the right hand side is the approximate CDF of $\psi$ bounded by $q_m$.  The summation is an approximate CDF as it is dependent on a kernel designed such that $K(\eta)$ is bounded by unity.  Finally, using the approximation of Eq.~\eqref{eq:CDF}, the KDE is related to the chance constraint of Eq.~\eqref{eq:CC_scalar} as follows:
\begin{equation}\label{eq:KDEeq}
P(\psi \geq q_m) \approx 1-\frac{1}{N} \sum_{j = 1}^N K(\eta_j). 
\end{equation}
Comparing Eq.~\eqref{eq:KDEeq} to Eq.~\eqref{eq:equalApp}, the approximate chance constraint resulting from Method 2 does not bound the original chance constraint like the approximation from Method 1.  As a result, the disadvantage of Method 2 is that the approximate chance constraint resulting from the application of KDEs can lie outside the boundary of the original chance constraint.  

\subsubsection{Comparison of the Two Methods}\label{sect:newm}
When comparing the two methods described in Sections~\ref{sect:firstm} and~\ref{sect:secondm}, the method of Section~\ref{sect:firstm} has the advantage of placing an upper bound on the chance constraint.  The functions used in this method have the requirement of Eq.~\eqref{eq:boundind}, leading to an approximation of the CDF of $\psi$ that is not bounded by unity.  It is noted, however, that KDEs are designed with functions that lead to an approximate CDF.  Combining this last fact with the wide variety of kernels that can be applied to a given problem, makes the method of Section~\ref{sect:secondm} an attractive option for approximating a chance constraint.  The main issue with the method of Section~\ref{sect:secondm} is the possibility of violating the bounds on the chance constraint.  The new method is an improvement on the methods of Sections~\ref{sect:firstm} and~\ref{sect:secondm} as it utilizes the necessary underlying functional form to ensure the upper bound on the chance constraint from the method of Section~\ref{sect:firstm}, while retaining the versatility and other characteristics of KDEs.  In the new method, MCMC samples of $\g{\xi}$ are used with KDEs that satisfy the necessary requirements to obtain a new nonlinear approximation of the chance constraint that upper bounds the chance constraint.   A nonlinear approximation stemming from Method 1 that both upper bounds the chance constraint and with the right choice of parameters is also a KDE is the Split-Bernstein approximation.  Consequently, to determine the necessary requirements for a KDE to upper bound the chance constraint, the Split-Bernstein approximation will be studied.  

\subsection{Split-Bernstein Approximation}\label{sect:SBkernel}
\begin{figure}[ht!]
  \centering
  \vspace*{0.25in}
{\includegraphics[height = 2.5in]{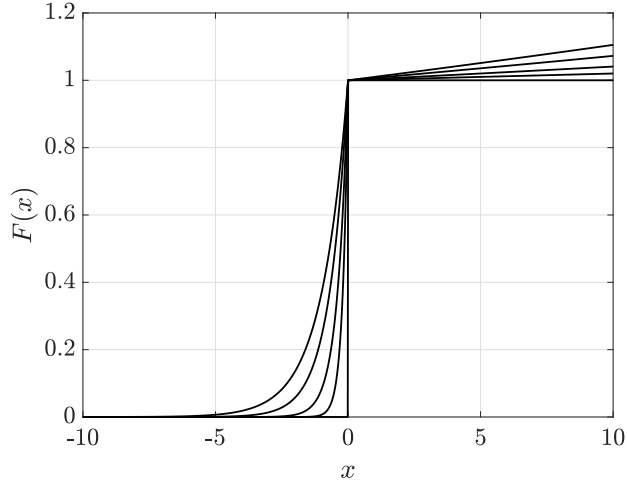}}
\caption{Convergence of the Split-Bernstein function as $\alpha_+ \rightarrow 0$ and $\alpha_- \rightarrow + \infty$.}
\label{fig:Kernel}
\end{figure}

In this section, the Split-Bernstein approximation and its KDE form are presented to define the necessary criteria for a KDE to upper bound the chance constraint.  The Split-Bernstein approximation from Ref.~\cite{Kumar1} follows from the inequality of Eq.~\eqref{eq:equalApp} as: 
\begin{equation}\label{eq:upperBou}
1-\epsilon_m \leq \bb{E}_{\g{\xi}} [1_{[q_m,+\infty)}(\alpha (\psi-q_m))] \leq \bb{E}_{\g{\xi}} [\Xi_{\alpha}(\alpha (\psi-q_m))], \ \forall \alpha>0,
\end{equation}
where the function $w(\cdot)$ from Eq.~\eqref{eq:equalApp} is replaced by the following Split-Bernstein function $\Xi_{\alpha}$: 
\begin{equation} \label{eq:piecewise}
\Xi_{\alpha} (\alpha ( \psi-q_m) ) = 
\begin{cases}
\exp(\alpha_+  (\psi-q_m)), \ \textrm{if} \ \psi-q_m  \geq 0, \\
 \exp(\alpha_- (\psi-q_m)), \ \textrm{if} \ \psi-q_m  < 0.
 \end{cases}
\end{equation}
Both $\alpha_+$ and $\alpha_-$ are parameters defined such that the piecewise function of Eq.~\eqref{eq:piecewise} converges to the indicator function when $\alpha_+ \rightarrow 0$ and $\alpha_- \rightarrow + \infty$.  This convergence can be seen in Fig.~\ref{fig:Kernel} where $\psi \geq x$ is an event with bound $x = 0$ and $F(x)$ represents the Split-Bernstein function.

The Split-Bernstein upper bound on the constraint from Eq.~\eqref{eq:upperBou} requires evaluation of the expectation, where computing the expectation requires integration of the unknown PDF of $\psi$.  Samples of $\psi(\g{\xi})$ are available as a function of the MCMC samples of $\g{\xi}$.  Based on the samples of $\g{\xi}$, an expression that converges to the expectation $\bb{E}_{\g{\xi}} [\Xi_{\alpha}(\cdot)]$ of Eq.~\eqref{eq:upperBou} was developed in Ref.~\cite{Kumar1}.
Using the resulting expression from Ref.~\cite{Kumar1}, the Split-Bernstein expectation of Eq.~\eqref{eq:upperBou} can be approximated as:
\begin{equation}\label{eq:SB_MCMC}
 \frac{1}{N} \sum_{j \in J_+} \exp (\alpha_+  (q_m-\psi_j)) +  \frac{1}{N} \sum_{j \in J_-} \exp (\alpha_-  (q_m-\psi_j)) \leq \epsilon_m.
\end{equation} 
For the summation of Eq.~\eqref{eq:SB_MCMC}, the event $\psi < q_m$ is indexed as $ J_+ \in \left( q_m-\psi_j  > 0 \right) $ and $ J_- \in  \left( q_m-\psi_j  \leq 0 \right)$.  The summation on the left hand side of Eq.~\eqref{eq:SB_MCMC} is the Split-Bernstein approximation.

\subsubsection{Split-Bernstein Approximation as a KDE}\label{sect:SB_KDE}
The Split-Bernstein approximation in Eq.~\eqref{eq:SB_MCMC} assigns a number that can be greater than unity for samples of $\psi$ in $ J_+ \in \left( q_m-\psi_j  > 0 \right) $.  As such, the Split-Bernstein approximation is not an approximate CDF in its general form.  KDEs employ kernels that are bounded by unity, thus ensuring that the right hand side of Eq.~\eqref{eq:CDF} derived from the general form of a KDE is an approximate CDF.  Choosing $\alpha_+ = 0$ and $\alpha_- = \alpha$ for the Split-Bernstein approximation ensures that the approximation cannot exceed unity.  Therefore, Eq.~\eqref{eq:SB_MCMC} can be formulated as the following approximate CDF:
\begin{equation}\label{eq:SB_alf}
\frac{1}{N} \sum_{j \in J_+} 1 +  \frac{1}{N} \sum_{j \in J_-} \exp (\alpha  (q_m-\psi_j)) \leq \epsilon_m,
\end{equation}
where the left-hand side of Eq.~\eqref{eq:SB_alf} is the Split-Bernstein approximate CDF.  Consequently, taking the derivative of the left hand side of Eq.~\eqref{eq:SB_alf} leads to the following approximate PDF of $\psi$:
\begin{equation}\label{eq:PDF_SB}
\hat{f}_{\psi}(q_m) = \frac{1}{N} \sum_{j \in J_+} 0 +  \frac{1}{N} \sum_{j \in J_-} \alpha \exp (\alpha  (q_m-\psi_j)).
\end{equation}
Suppose now that $\alpha$ is defined as
\begin{equation}\label{eq:alpdes}
\alpha = \frac{1}{h}.
\end{equation}
Using $\alpha$ from Eq.~\eqref{eq:alpdes}, Eq.~\eqref{eq:PDF_SB} transforms to the form of Eq.~\eqref{eq:KDE}, indicating that the Split-Bernstein approximation is a potential KDE.  The associated candidate kernel $k_{SB}(\cdot)$ from Eq.~\eqref{eq:PDF_SB} for the Split-Bernstein approximation is
 \begin{equation}\label{eq:SB_KDE}
k_{SB}(\eta) = 
\begin{cases}
0, \  \eta > 0 \\
\exp(\eta), \ \eta \leq 0.
\end{cases}
\end{equation}
The inequality cases of Eq.~\eqref{eq:SB_KDE} replace the inequalities of Eq.~\eqref{eq:piecewise} by taking
\begin{equation}\label{eq:etaproof}
\eta > 0 \  \Rightarrow \ \frac{q_m-\psi}{h} > 0 \ \Rightarrow  \ q_m -\psi> 0.
\end{equation}
It remains to show that the Split-Bernstein candidate kernel satisfies the requirements of a kernel.  The requirements of a kernel vary in the literature~\cite{Hwang,Mittal}, so three main criteria were chosen by the authors as they encompass the most common requirements. These criteria lead to the following proposition.

{\flushleft \m{Proposition 1.}} The Split-Bernstein approximation with $\alpha_+ = 0$ leads to the kernel of Eq.~\eqref{eq:SB_KDE} that satisfies the following criteria:
\begin{enumerate}
\item $\displaystyle  k(\eta) \geq 0$, 
\item $\int_{-\infty}^{+\infty} k(x) dx = 1$, 
\item $0 < \int_{-\infty}^{+\infty} x^2 k(x)dx < +\infty$.
\end{enumerate}

{\flushleft \m{Proof of Proposition 1.}}  For Case (1), the Split-Bernstein candidate kernel is zero for $\eta > 0$ and always greater than zero for $\eta \leq 0$.  Showing that Case (2) is satisfied is straightforward:
\begin{equation}\label{eq:proof2}
\int_{-\infty}^{+\infty} k_{SB}(x) dx \Rightarrow \int_{-\infty}^0 \exp(x) dx + \int_0^{+\infty} 0 dx \Rightarrow  \exp(x) \bigg\rvert_{-\infty}^0 + 0 = 1.
\end{equation}
Likewise, Case (3) is proven as:
\begin{equation}\label{eq:proof3}
\int_{-\infty}^{+\infty} x^2 k_{SB}(x)dx \Rightarrow \int_{-\infty}^0 x^2 \exp(x) dx + \int_0^{+\infty} 0 dx = 2.
\end{equation}
\qed

The Split-Bernstein candidate kernel $k_{SB}$ is thus shown to be a kernel.  As a result of this fact, the Split-Bernstein approximation is a KDE for $\alpha_+ = 0$ and $\alpha_- = \alpha$.  Due to the use of MCMC samples and the form of the Split-Bernstein kernel, the Split-Bernstein approximation retains an upper bound on the chance constraint.  Thus in Section~\ref{sect:BiasKDE}, the Split-Bernstein kernel is analyzed to determine the correct criteria for a kernel that, when combined with the use of MCMC samples, ensures that the associated KDE retains the same upper bound on the chance constraint as Method 1 from Section~\ref{sect:firstm}.

\subsection{New Method: Biased Kernel Density Estimators}\label{sect:BiasKDE}
The approximate probability form derived from a KDE is dependent on the choice of kernel $k(\eta)$.  A kernel is placed at each sample and the approximate probability is the average of the integrated sum of these kernels evaluated at the weighted sample distance from some bound $q_m$ (see Eq.~\eqref{eq:KDEeq}).  Kernels are often designed such that their mean is at the same location as the sample.  Having the mean shift to the right or left of the sample is referred to as ``biasing'' in this paper.  As the approximate probability is the result of an average that is a function of the kernel, changing the location of the kernel mean will cause the approximate probability to increase or decrease.  Biasing the kernel will be shown in the following section to be one component of reformulating a KDE so that it will provide an upper bound on the chance constraint.  

\subsubsection{Bias of the Split-Bernstein Kernel}\label{sect:SBbias}
\begin{figure}[ht!]
  \centering
  \vspace*{0.25in}
{\includegraphics[height = 2.5in]{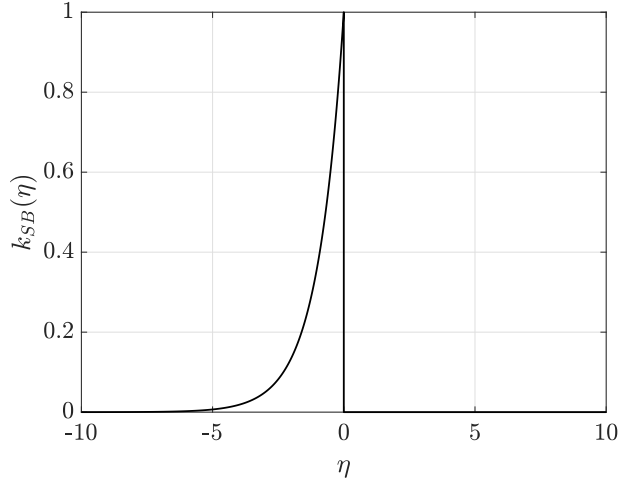}}
\caption{Split-Bernstein kernel.}
\label{fig:KernelSB}
\end{figure}

In this section, the Split-Bernstein kernel is studied to determine whether or not it is a biased kernel.  Figure~\ref{fig:KernelSB} shows the Split-Bernstein kernel for a sample centered at zero.  It is seen that the Split-Bernstein kernel has a sail like shape from $-\infty$ to zero, and is zero otherwise.  This shape implies that the mean would be shifted to the left of the sample, indicating that the Split-Bernstein kernel is biased.  \m{Proposition 2} states that the Split-Bernstein kernel is a biased kernel, and this proposition is subsequently proven.  

%
{\flushleft \m{Proposition 2.}} The Split-Bernstein kernel is left biased with a mean $\mu = -1$, therefore the Split-Bernstein kernel must be shifted by the bandwidth $h$ to have $\mu = 0$.  

{\flushleft \m{Proof of Proposition 2.}}  First, using the definition of the expectation:
\begin{equation}\label{eq:biasproof}
\bb{E}[k_{SB}(x)] = \int_{-\infty}^{0} x \exp(x) dx = -1, 
\end{equation}
proving that $\mu = -1$.  Due to the design of the Split-Bernstein kernel, to make the kernel unbiased the mean must be shifted by some constant $c$:
\begin{equation}\label{eq:solveforc}
\int_{-\infty}^{c} x \exp(x-c) dx = 0, \ \Rightarrow \int_{-\infty}^{c} x \exp(x-c) dx = c-1 = 0, \ \Rightarrow c = 1. 
\end{equation}
Because the integral of Eq.~\eqref{eq:kernCDF} was defined using a bound transformation in Eq.~\eqref{eq:CDFactual}, the constant $c$ has the following relation:
\begin{equation}\label{eq:solvebias}
c = \frac{q_m - \psi_j}{h}, \ \Rightarrow \frac{q_m - \psi_j}{h} = 1, \ \Rightarrow q_m - \psi_j = h,
\end{equation}
proving that the Split-Bernstein kernel must be shifted by the bandwidth $h$ to be unbiased.
\qed

In order to demonstrate the importance of the shift by $h$ from \m{Proposition 2}, consider biasing the Split-Bernstein kernel using a mean shift of $\delta$ where 
\begin{equation}\label{eq:definedelt}
0 < \delta \leq  \frac{1}{\alpha}, \ \delta \in \bb{R}.
\end{equation}
It is noted that the relation of $h$ and $\alpha$ from Eq.~\eqref{eq:alpdes} has been applied in Eq.~\eqref{eq:definedelt}.  The shift $\delta$ results in the following biased Split-Bernstein kernel function $K_{B,SB}(\cdot)$:
\begin{equation} \label{eq:piecewiseNew}
K_{B,SB} (\alpha  (\psi-q_m) ) = 
\begin{cases}
1, \ \textrm{if} \ \psi-q_m  \geq \delta, \\
\exp(\alpha (\psi-q_m-\delta)), \ \textrm{if} \  \delta > \psi-q_m  \geq  0, \\
 \exp(\alpha (\psi-q_m-\delta)), \ \textrm{if} \ \psi-q_m  < 0.
 \end{cases}
\end{equation}
Recalling from Eq.~\eqref{eq:boundind}, the Split-Bernstein function has the property: 
\begin{equation}\label{eq:functinequal}
1_{[q_m,+\infty)}(\alpha (\psi-q_m)) \leq \Xi_{\alpha}(\alpha (\psi-q_m)), \ \ \forall \psi.
\end{equation}
The relation between the Split-Bernstein function, the function $K_{SB}(\cdot)$, and the indicator function can be obtained by application of Eqs.~\eqref{eq:kernCDF} and~\eqref{eq:SB_KDE} with $\alpha_+ = 0$ and $\alpha_- = \alpha$:
\begin{equation}\label{eq:SBequaK}
1_{[q_m,+\infty)}(\alpha (\psi-q_m)) \leq \Xi_{\alpha}(\alpha (\psi-q_m)) = K_{SB} (\alpha (\psi-q_m)).
\end{equation}
For comparison of the indicator function with the new function $K_{B,SB} (\cdot)$ in Eq.~\eqref{eq:piecewiseNew}, the indicator function is defined as:
\begin{equation} \label{eq:indicator}
1_{[q_m,+\infty)} (\alpha  (\psi-q_m) ) = 
\begin{cases}
1, \ \textrm{if} \ \psi  \geq q_m, \\
0, \ \textrm{if} \ \psi  < q_m.
 \end{cases}
\end{equation}
Comparing the function $K_{B,SB} (\cdot)$ in Eq.~\eqref{eq:piecewiseNew} with the indicator function in Eq.~\eqref{eq:indicator}, the inequality of Eq.~\eqref{eq:SBequaK} no longer applies as a portion of the function $K_{B,SB} (\cdot)$ lies below that of the indicator function.  As a result, the conservative bound of Eq.~\eqref{eq:equalApp} can no longer be guaranteed.  Thus a mean shift that is equal to or smaller than zero is required to ensure that the Split-Bernstein KDE is an upper bound on the chance constraint of Eq.~\eqref{eq:equalApp}.  Additionally, if a mean shift is applied to the kernel, it is necessary that it be a function of the bandwidth $h$ to ensure that the kernel function $K(\cdot)$ will still converge to the indicator function as the bandwidth goes to zero.

\subsubsection{Biased Form of General KDEs}\label{sect:GenBiasKDE}
In the previous section it was determined that the bias of the Split-Bernstein kernel, when combined with MCMC sampling, ensures the upper bound on the chance constraint from Eq.~\eqref{eq:equalApp}.  This result can be extended to other KDEs that satisfy the requirements of the following \m{Theorem 1}. 

{\flushleft \m{Theorem 1.}}  Suppose there exists a biased kernel $k_B(\nu)$ defined such that $k_B(\nu)$ has an integrated function $K_B(\nu)$ that satisfies the inequality 
\begin{align}\label{eq:geninequality}
1_{[y,+\infty)} \left( \nu \right) \leq K_B(\nu), \ \forall \nu, ~ \textrm{with}~h >0, \\
\nu = \frac{Y-y}{h}
\end{align} 
for some event $Y \geq y$ in probability space $P(\cdot)$ where $Y$ is a random variable supported on set $\Omega \subseteq \bb{R}$.  The biased kernel $k_B(\nu)$ has the following relation to the unbiased kernel $k(\nu)$:
\begin{equation}
\bb{E} [ K_{B} (\nu)] = \bb{E} [ K (\nu)] \pm \frac{B(h)}{h}, \ B(h) > 0,
\end{equation}
where the bias is set as $B(h)$.  If enough MCMC samples of the random variable $Y_j$ are available to satisfy the law of large numbers, 
the following inequality holds: 
\begin{equation}\label{eq:kernelCDFbound}
\frac{1}{N} \sum_{j = 1}^N \int_{-\infty}^{\eta_j}k_B(x_j) d x_j \leq P(Y<y) \leq \epsilon_m,
\end{equation}
where $\epsilon_m\in \bb{R}$ and $\eta$ was defined in Eq.~\eqref{eq:etadef}.

{\flushleft \m{Proof of Theorem 1.} } Suppose that a PDF $f_{Y}(\cdot)$ exists and has bounded first and second moments with a CDF $F_{Y}(\cdot)$ satisfying:
\begin{equation}\label{eq:needed}
1_{[y,+\infty)} \left( \nu \right) \leq F_{Y}(\nu), \forall \nu,
\end{equation}
where MCMC samples $Y_j$ of random variable $Y$ are available.  Then the sample expectation $\hat{\bb{E}}$ is defined as
\begin{equation}\label{eq:sampleExp}
\hat{\bb{E}} [F_{Y}(\nu)] = \frac{1}{N} \sum_{j = 1}^N F_{Y} (\nu_j).
\end{equation}
Now define the following nonempty compact set:
\begin{equation}\label{eq:setsforproof}
1-\epsilon_m \leq \bb{E} [ F_{Y} (\nu)], \\
\end{equation}
where $\bb{E} [ F_{Y} (\nu)]$ exists for $h > 0$.  
Assuming that enough MCMC samples $Y_j$ of the random variable are available, leads to the following relation:
\begin{equation}\label{eq:inequaltosat}
\frac{1}{N} \sum\limits_{j = 1}^N F_{Y} (\nu_j)  \ \Rightarrow \  \bb{E} [ F_{Y} (\nu)] \geq 1-\epsilon_m. \\
\end{equation}
Independent of the result of Eq.~\eqref{eq:inequaltosat}, it is a known property of probability that
\begin{equation}\label{eq:probind}
P(Y \geq y) = \bb{E} \left[ 1_{[y,+\infty]} \left( Y-y) \right) \right] = \bb{E} \left[ 1_{[y,+\infty)} \left( \frac{Y-y}{h} \right) \right], \ h> 0,
\end{equation}
where scaling the indicator function by $h$ does not affect the relation.  Additionally, due to the relation of Eq.~\eqref{eq:needed} the following property of probability can be applied to Eq.~\eqref{eq:probind}:
\begin{equation}\label{eq:prob}
 P(Y \geq y) = \bb{E} \left[ 1_{[y,+\infty)} \left( \nu \right) \right] \leq \bb{E} [F_{Y} (\nu)].
\end{equation}
Substituting the relation from Eq.~\eqref{eq:inequaltosat} into Eq.~\eqref{eq:prob} leads to the following inequality:
\begin{equation}\label{eq:probcdf}
1-\epsilon_m \leq P(Y \geq y) \leq \frac{1}{N} \sum_{j = 1}^N F_{Y} (\nu_j).
\end{equation}
The form of the PDF $f_{Y}(\cdot)$ is such that it satisfies the requirements for a kernel from \m{Proposition 1}.  Therefore 
\begin{equation}\label{eq:kerntoprob}
F_{Y}(\nu) = K_B(\nu).
\end{equation}
Substituting the relation of Eq.~\eqref{eq:kerntoprob} into Eq.~\eqref{eq:probcdf} leads to the following relation:
\begin{equation}\label{eq:probcdfK}
1-\epsilon_m \leq P(Y \geq y) \leq \frac{1}{N} \sum_{j = 1}^N K_B(\nu_j).
\end{equation}
Changing the probability form to its complement in Eq.~\eqref{eq:probcdfK} and substituting $\eta$ from Eq.~\eqref{eq:etadef} leads to the following relation:
\begin{equation}\label{eq:compprob}
\frac{1}{N} \sum_{j = 1}^N K_B \left( \eta_j \right) \leq P(Y<y) \leq \epsilon_m.
\end{equation}
Applying the relation from Eq.~\eqref{eq:kernCDF} to Eq.~\eqref{eq:compprob} results in the following inequality:
\begin{equation}\label{eq:proofresult}
\frac{1}{N} \sum_{j = 1}^N \int_{-\infty}^{\eta_j}k_B(x_j) d x_j \leq P(Y<y) \leq \epsilon_m.
\end{equation}
This concludes the proof.\qed

Thus the first requirement for a kernel to result in a KDE that provides an upper bound on the chance constraint is that the kernel is biased by a bandwidth dependent amount that is enough to ensure the inequality of Eq.~\eqref{eq:needed}.  The second requirement is that the sampling method applied ensures convergence of the sample expectation to the expectation that upper bounds the chance constraint from Eq.~\eqref{eq:equalApp}.  The new method of biased KDEs applies kernels that satisfy the requirements of \m{Theorem 1}, as well as MCMC sampling to obtain the nonlinear approximation of the chance constraint from Eq.~\eqref{eq:kernelCDFbound}.  This nonlinear approximation provides an upper bound on the chance constraint as well as being an approximate CDF of the random variable, thus retaining the advantages of Methods 1 and 2 from Sections~\ref{sect:firstm} and~\ref{sect:secondm}, respectively, in a generally applicable manner.  

\subsection{Biasing the Epanechnikov and Gaussian Kernels}\label{sect:commonKDE}
Section~\ref{sect:BiasKDE} established the requirement for a kernel to result in a KDE that provides an upper bound on the chance constraint.  In this section, this requirement will be applied to the Epanechnikov and Gaussian kernels as these are two common kernels.  

\subsubsection{Epanechnikov Kernel}\label{sect:epanechbias}
The Epanechnikov kernel is considered as it is often referenced as the optimal kernel due to its convergence properties~\cite{Hodges}.  Letting $B$ represent a bias added to the Epanechnikov kernel from Ref.~\cite{Epanech1}, leads to the following representation of the biased kernel function $K_{B,E}$:   
\begin{equation} \label{eq:epanechkernel}
K_{B,E} (\nu_j ) = 
\begin{cases}
1, \ \textrm{if} \ \nu_j  \geq 1 - \frac{B}{h_E}, \\
\frac{1}{2}+\frac{3}{4} \left( \nu_j+\frac{B}{h_E} \right) -\frac{1}{4} \left( \nu_j +\frac{B}{h_E} \right) ^3, \ \textrm{if} \ -1-\frac{B}{h_E} < \nu_j  <  1-\frac{B}{h_E},  \\
 0, \ \textrm{if} \ \nu_j  \leq -1-\frac{B}{h_E},
 \end{cases}
\end{equation}
where $h_E$ is the bandwidth of the Epanechnikov kernel and $\nu_j$ is defined as
\begin{equation}
\nu_j = \frac{\psi_j - q_m}{h}.
\end{equation}
 Setting $B$ equal to the bandwidth $h_E$ ensures that the kernel function $K_{B,E}$ remains greater than or equal to the indicator function everywhere and so satisfies the inequality of Eq.~\eqref{eq:geninequality}.

\subsubsection{Gaussian Kernel}\label{sect:gausbias}
\begin{figure}[ht!]
  \centering
  \vspace*{0.25in}
{\includegraphics[height = 2.5in]{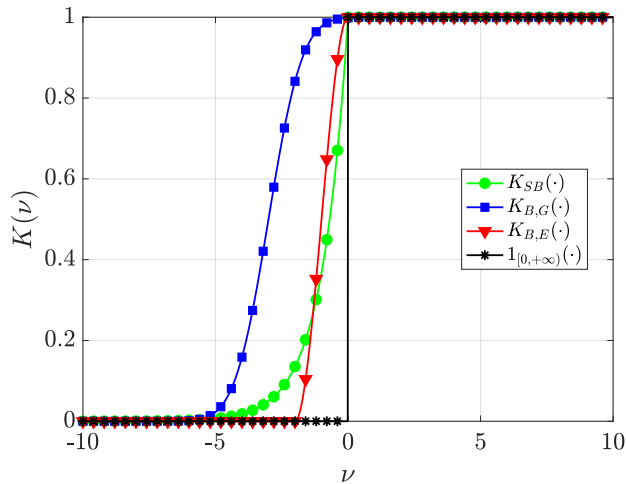}}
\caption{Biased kernel function $K_B(\cdot)$ and the indicator function.}
\label{fig:Kernel_Compare}
\end{figure}

The Gaussian kernel is a popular kernel as it is completely smooth and is fairly straightforward to implement numerically.  The design of the biased Gaussian kernel is
\begin{equation}\label{eq:gausskern}
K_{B,G} (\nu_j) = \frac{1}{2} \left( 1+ \erf \left( \frac{ h_G \nu_j+B}{h_G \sqrt{2}} \right) \right), \ \forall \nu_j,
\end{equation}
where $h_G$ is both the bandwidth and the standard deviation and $B$ represents the added bias.  Solving for the bias that will result in $K_{B,G} = 1$: 
\begin{equation}\label{eq:findGbias}
K_{B,G} (\nu_j) = 1 \ \Rightarrow \ \erf \left( \frac{ h_G \nu_j + B}{h_G \sqrt{2}} \right) = 1 \ \Rightarrow \ \frac{ h_G \nu_j + B}{h_G \sqrt{2}} = \erf^{-1}(1) \ \Rightarrow \ B = \infty.
\end{equation}
As a result, the Gaussian kernel cannot be biased enough to ensure the inequality of Eq.~\eqref{eq:geninequality}.  It is possible to ensure that at least $99 \% $ of the Gaussian kernel function $K_{B,G}$ satisfies the requirement of Eq.~\eqref{eq:geninequality} by assigning a bias of $3 h_G$.  Even though the Gaussian kernel does not satisfy the requirements of \m{Theorem 1}, it could still provide a useful comparison to other less smooth kernels in numerical application.  The three kernel functions $K (\cdot)$ that have been discussed thus far are shown for comparison in Fig.~\ref{fig:Kernel_Compare} along with the indicator function.   It is noted that the bias for the Split-Bernstein function $K_{SB}$ and the indicator function are both zero due to the design of the functions.  

The new method of biased KDEs that leads to a nonlinear approximation of the chance constraint that provides an upper bound on the chance constraint has been derived.  This method can be applied to transform chance constrained optimization problems to deterministic optimization problems, particularly chance constrained optimal control problems.    

\section{Chance Constrained Optimal Control \label{sect:bolza}}
Consider the following continuous time CCOCP.  Determine the state
$\m{y}(\tau)\in\bb{R}^{n_y}$ and the control $\m{u}(\tau)\in\bb{R}^{n_u}$ on the domain $\tau \in [-1, +1]$,  the initial time, $t_0$, and the terminal time $t_f$ that minimize the cost functional

\begin{subequations}
\begin{equation}\label{bolza-cost-s}
  \C{J} =\C{M}(\m{y}(-1),t_0,\m{y}(+1),t_f)  + \frac{t_f-t_0}{2}\int_{-1}^{+1} \C{L}(\m{y}(\tau),\m{u}(\tau), t(\tau, t_0, t_f))\, d\tau,
\end{equation}
subject to the dynamic constraints
\begin{equation}\label{bolza-dyn-s}
  \frac{d\m{y}}{d\tau} -
\frac{t_f-t_0}{2}\m{a}(\m{y}(\tau),\m{u}(\tau), t(\tau, t_0, t_f) )=\m{0}, 
\end{equation}
the inequality path constraints
\begin{equation}\label{bolza-path-s}
\m{c}_{\min} \leq \m{c}(\m{y}(\tau),\m{u}(\tau), t(\tau, t_0, t_f) )\leq \m{c}_{\max},
\end{equation}
the boundary conditions
\begin{equation}\label{bolza-bc-s}
  \m{b}_{\min} \leq \m{b}(\m{y}(-1),t_0,\m{y}(+1),t_f) \leq \m{b}_{\max},
\end{equation}
and the chance-constraints
\begin{equation}\label{bolza-pathcc-s}
P( \m{F} (\m{y}(\tau),\m{u}(\tau),t(\tau, t_0, t_f);\g{\xi}) \geq \m{q}) \geq 1- \epsilon.
\end{equation}
\end{subequations}

It is noted that the time interval $\tau\in[-1,+1]$ can be transformed
to the time interval $t\in[t_0,t_f]$ via the affine transformation
\begin{equation}\label{tau-to-t}
  t \equiv t(\tau,t_0,t_f) = \frac{t_f-t_0}{2}\tau + \frac{t_f+t_0}{2}.  
\end{equation}
where $\m{y}$ is the state defined on the feasible set $\m{Y} \subset \bb{R}^n$. The control $\m{u}$ is defined on the feasible set $\m{U} \subset \bb{R}^m$.   The parameters and variables related to the chance constraint were described in Section~\ref{sect:SBkernel}.  It should be noted that a chance constraint can be used to represent both path, dynamic and event probabilistic constraints.

Before applying the method of biased KDEs to CCOCPs to transform the chance constraints to nonlinear constraint approximations, the continuous time CCOCP must be transformed to a form that can be solved using numerical methods.  For application with numerical methods, the CCOCP is discretized on the domain $\tau\in[-1,+1]$ which is partitioned into a {\em mesh} consisting of $K$  {\em mesh intervals} $\C{S}_k=[T_{k-1},T_k],\; k=1,\ldots,K$, where $-1 = T_0 < T_1 < \ldots < T_K = +1$.  The mesh intervals have the property that $\displaystyle \cup_{k=1}^{K} \C{S}_k=[-1,+1]$.  Let $\m{y}^{(k)}(\tau)$ and $\m{u}^{(k)}(\tau)$ be the state and control in $\C{S}_k$.  Using the transformation given in Eq.~\eqref{tau-to-t}, the chance constrained optimal control problem of Eqs.~\eqref{bolza-cost-s}-\eqref{bolza-pathcc-s} can then be rewritten as follows. Minimize the cost functional

\begin{subequations}
\begin{equation}\label{bolza-cost-segmented}
    \C{J} = \C{M}(\m{y}^{(1)}(-1),t_0,\m{y}^{(K)}(+1),t_f)  + \frac{t_f-t_0}{2}\sum_{k=1}^K \int_{T_{k-1}}^{T_k} \C{L}(\m{y}^{(k)}(\tau),\m{u}^{(k)}(\tau),t)\, d\tau,
\end{equation}
subject to the dynamic constraints
\begin{equation}\label{bolza-dyn-segmented}
  \displaystyle\frac{d\m{y}^{(k)}(\tau)}{d\tau} - \frac{t_f-t_0}{2}\m{a}(\m{y}^{(k)}(\tau),\m{u}^{(k)}(\tau), t)=\m{0}, \quad  (k=1,\ldots,K),
\end{equation}
the path constraints
\begin{equation}\label{bolza-path-segmented}
  \m{c}_{\min} \leq \m{c}(\m{y}^{(k)}(\tau),\m{u}^{(k)}(\tau), t) \leq \m{c}_{\max},\quad  (k=1,\ldots,K),
\end{equation}
the boundary conditions
\begin{equation}\label{bolza-bc-segmented}
\m{b}_{\min} \leq \m{b}(\m{y}^{(1)}(-1),t_0,\m{y}^{(K)}(+1),t_f) \leq  \m{b}_{\max}, 
\end{equation}
and the chance-constraints
\begin{equation}\label{bolza-pathcc-segmented}
P(\m{F} (\m{y}^{(k)}(\tau),\m{u}^{(k)}(\tau),t;\g{\xi}) \geq \m{q}) \geq 1- \epsilon.
\end{equation}
\end{subequations}
Because the state must be continuous at each interior mesh point, it
is required that the condition $\m{y}(T_k^{-})=\m{y}(T_k^{+}),\;(k=1,\ldots,K-1)$ be
satisfied at the interior mesh points $(T_1,\ldots,T_{K-1})$.

\section{Legendre-Gauss-Radau Collocation\label{sect:LGR}}
The form of discretization that will be applied to the CCOCP in Section~\ref{sect:bolza} is collocation at
Legendre-Gauss-Radau (LGR) points~\cite{Garg1,Garg2,Patterson2015}.  In the LGR 
collocation method, the state of the continuous-time CCOCP is approximated in $\C{S}_k,\;k\in[1,\ldots,K]$, as 
\begin{equation}\label{state-approximation-LGR}
\begin{split}
\m{y}^{(k)}(\tau)  \approx \m{Y}^{(k)}(\tau) & = \sum_{j=1}^{N_k+1}
\m{Y}_{j}^{(k)} \ell_{j}^{(k)}(\tau),\\ \ell_{j}^{(k)}(\tau) & = \prod_{\stackrel{l=1}{l\neq j}}^{N_k+1}\frac{\tau-\tau_{l}^{(k)}}{\tau_{j}^{(k)}-\tau_{l}^{(k)}}, 
\end{split}
\end{equation}  
where $\tau\in[-1,+1]$, $\ell_{j}^{(k)}(\tau),$ $j=1,\ldots,N_k+1$, is a
basis of Lagrange polynomials,
$\left(\tau_1^{(k)},\ldots,\tau_{N_k}^{(k)}\right)$ are the 
Legendre-Gauss-Radau (LGR)~\cite{Garg1} collocation points
in $\C{S}_k =$ $[T_{k-1},T_k)$, and 
$\tau_{N_k+1}^{(k)}=T_k$ is a noncollocated point.  Differentiating
$\m{Y}^{(k)}(\tau)$ in Eq.~(\ref{state-approximation-LGR}) with
respect to $\tau$ gives
\begin{equation}\label{diff-state-approximation-LGR}
  \frac{d\m{Y}^{(k)}(\tau)}{d\tau} = \sum_{j=1}^{N_k+1}\m{Y}_{j}^{(k)}\frac{d\ell_j^{(k)}(\tau)}{d\tau}.
\end{equation}
Defining $t_i^{(k)}=t(\tau_i^{(k)},t_0,t_f)$ using
Eq.~\eqref{tau-to-t}, the dynamics are then approximated at the $N_k$
LGR points in mesh interval $k\in[1,\ldots,K]$ as
\begin{equation}\label{collocation-LGR}
    \sum_{j=1}^{N_k+1}D_{ij}^{(k)} \m{Y}_j^{(k)} - \frac{t_f-t_0}{2}\m{a}(\m{Y}_i^{(k)},\m{U}_i^{(k)},t_i^{(k)})=\m{0}, \  (i=1,\ldots,N_k),
\end{equation}
where $D_{ij}^{(k)} = d\ell_j^{(k)}(\tau_i^{(k)})/d\tau,\;(i=1,\ldots,N_k),\;(j=1,\ldots,N_k+1)$ are the elements of the $N_k\times (N_k+1)$ {\em Legendre-Gauss-Radau differentiation matrix}~\cite{Garg1} in mesh interval
$\C{S}_k,\;k\in[1,\ldots,K]$.  The LGR discretization then leads to
the following programming problem.  Minimize 
\begin{equation}\label{cost-LGR}
    \C{J}  \approx \C{M}(\m{Y}_{1}^{(1)},t_0,\m{Y}_{N_K+1}^{(K)},t_f) +   \sum_{k=1}^{K} \sum_{j=1}^{N_k}  \frac{t_f-t_0}{2} 
  w_{j}^{(k)} \C{L}(\m{Y}_{j}^{(k)},\m{U}_{j}^{(k)},t_j^{(k)}),
\end{equation}
subject to the collocation constraints of Eq.~(\ref{collocation-LGR})
and the constraints
\begin{gather}\label{eq:differential-collocation-conditions-LGR}
  \m{c}_{\min} \leq \m{c}(\m{Y}_{i}^{(k)},\m{U}_{i}^{(k)},t_i^{(k)}) \leq \m{c}_{\max},\; (i=1,\ldots,N_k),\\
 \m{b}_{\min} \leq \m{b}(\m{Y}_{1}^{(1)},t_0,\m{Y}_{N_K+1}^{(K)},t_f)  \leq \m{b}_{\max},  \\
 P( \m{F} (\m{Y}_{i}^{(k)},\m{U}_{i}^{(k)},t_i^{(k)};\g{\xi}) \geq \m{q}) \geq 1- \epsilon ,\; (i=1,\ldots,N_k) , \label{cc-constraint}\\
\m{Y}_{N_k+1}^{(k)} = \m{Y}_1^{(k+1)} , \quad (k=1,\ldots,K-1),  \label{continuity-constraint}
\end{gather}
where $N = \sum_{k=1}^{K} N_k$ is the total number of LGR points and 
Eq.~\eqref{continuity-constraint} is the continuity condition on
the state and is enforced at the interior mesh points
$(T_1,\ldots,T_{K-1})$ by treating $\m{Y}_{N_k+1}^{(k)}$ and
$\m{Y}_1^{(k+1)}$ as the same variable in the programming problem.  

As the CCOCP is now in a programming problem form, the method of biased KDEs developed in Section~\ref{sect:CCnonlinear} can be applied to transform the chance constraint to a deterministic constraint.  As such, the chance constraint of Eq.~\eqref{cc-constraint} is transformed to the following deterministic constraint:
\begin{equation}\label{eq:LGRdetKDE}
\frac{1}{N} \sum_{j = 1}^N K_B \left( \eta_l \right) \leq \epsilon_m,
\end{equation}
The variable $\eta_l$ is defined as
\begin{equation}
\eta_l = \frac{q_m-F((\m{Y}_{i}^{(k)},\m{U}_{i}^{(k)},t_i^{(k)};\g{\xi}_l)}{h},
\end{equation}
where $l$ is the number of samples of the random variable $\g{\xi}$.  With the transformation of Eq.~\eqref{cc-constraint} to the deterministic constraint of Eq.~\eqref{eq:LGRdetKDE}, the programming problem is now a nonlinear programming problem (NLP) that can be solved using available software.

\section{Biased Kernel Method for Chance Constraints}\label{sect:method}
In this section, a computational framework to solve CCOCPs in their deterministic optimal control problem form is described that combines the theory of Sections~\ref{sect:CCnonlinear},~\ref{sect:bolza} and~\ref{sect:LGR}.  Such a framework can be computationally expensive due to the evaluation of the chance constraint approximation at each collocation point for 50,000+ MCMC samples.  For certain CCOCPs, a larger bandwidth can first be applied to obtain an over-smooth constraint approximation that can lead to faster convergence of the NLP.  When the NLP solver is in the neighborhood of the solution, the bandwidth can then be decreased to reach a less conservative solution.  The computational implementation of Sections~\ref{sect:CCnonlinear},~\ref{sect:bolza}, and~\ref{sect:LGR} is presented below as the Method of Biased Kernels for Chance Constrained Optimal Control.

\begin{shadedframe}
\vspace{-10pt}
\begin{center}
 \shadowbox{\bf Method of Biased Kernels for Chance Constrained Optimal Control}
\end{center}
\begin{enumerate}[{\bf Step 1:}]
\item Reformulate chance constraint to transform CCOCP to determinstic optimal control problem using the method of Section~\ref{sect:BiasKDE}.\label{step:KDE}
\item Obtain MCMC samples of random variables.\label{step:MCMC}
\item Choose a bandwidth $h$.\label{step:bandwidth}
\begin{enumerate}[{\bf (a):}]
\item Determine a trial $h_i$ value.\label{step:firstb}
\item Run optimal control problem through optimal control software for up to four mesh refinement iterations.\label{step:refine}
\item If NLP has converged and error tolerance for mesh is within $\epsilon$, set $h = h_i$, otherwise choose $h_{i+1} > h_{i}$ and return to {\bf~\ref{step:refine}}.
\end{enumerate}
\item Set bandwidth as $h > h_i$ if error tolerance decreases to $\epsilon$ in less iterations, based on {\bf Step~\ref{step:bandwidth}}, otherwise continue to {\bf Step~\ref{step:runfull}}.
\item Run optimal control problem through optimal control software with mesh refinement and possible switch of $h$ to $h_i$. \label{step:runfull}
\item Return to {\bf Step~\ref{step:runfull}} to repeat $20+$ times to ensure of NLP convergence for each run.
\end{enumerate}
\end{shadedframe}

\section{Application of Method}\label{sect:examples}
In this section, the computational framework described in Section~\ref{sect:method} is applied to an example.  The example is a chance constrained soft lunar landing optimal control problem from Ref.~\cite{Kumar1} that is a modification of the deterministic soft lunar landing optimal control problem from Ref.~\cite{Meditch1}.  The example is used to compare the application of the method of biased KDEs with three different kernels to results already available in the literature~\cite{Kumar1} for a CCOCP with a control dependent path chance constraint and an event chance constraint.  

The three kernels chosen for comparison were the Split-Bernstein, Epanechnikov, and Gaussian kernels.  Each kernel was chosen based on different properties of the integrated function $K(\cdot)$.  The Split-Bernstein function $K_{SB}(\cdot)$ tightly bounds the indicator function due to its left biased exponential form.  The Epanechnikov function $K_E(\cdot)$ does not provide as tight a bound, but it is a smoother function than the Split-Bernstein function $K_{SB}$ and was proven in Section~\ref{sect:GenBiasKDE} to satisfy the requirement of \m{Theorem 1} with the correct choice of bias.  The Epanechnikov function $K_E$ is a piecewise function like the Split-Bernstein function $K_{SB}$.  Thus the application of either kernel to available optimal control software could result in an increase in computational expense as opposed to a completely smooth kernel.  As a result, the Gaussian kernel was also included for comparison due to the smoothness of the function $K_G$, despite not having a bias that satisfies of the requirement of \m{Theorem 1} as discussed in Section~\ref{sect:BiasKDE}.  

The example was solved in MATLAB$^{\textrm{\textregistered}}$ version R2018a (build 9.4.0.813654) using the optimal control software package $\mathbb{GPOPS-II}$~\cite{Patterson2014} together with the NLP solver SNOPT~\cite{Gill1,Gill2}.  $\mathbb{GPOPS-II}$ implements an hp-adaptive LGR collocation method~\cite{Garg1,Garg2,Garg3,Patterson2015,Liu2015,Liu2018,Darby2,Darby3,Francolin2014a}.  Derivative approximations required by SNOPT were obtained using central finite-differencing~\cite{Patterson2012}.  Additionally, the NLP solver tolerance and maximum number of iterations were set to $10^{-6}$ and $500$, respectively.  Mesh refinement was performed using the method of Ref.~\cite{Liu2018} with a mesh refinement accuracy tolerance of $10^{-6}$ and decay rate of $0.5$ (as required by the method of Ref.~\cite{Liu2018}).  The initial mesh consisted of $10$ mesh intervals with four collocation points each.  A Hamiltonian Monte Carlo (HMC) method was chosen to obtain the $50,000$ MCMC samples as described in Neal~\cite{MCMCMethods,Neal2} that were generated for each run.  This number was used to ensure each set of MCMC samples generated per run was a convergent set. The bandwidths were determined using the MATLAB$^{\textrm{\textregistered}}$ function $ksdensity$ that employs methods from Refs.~\cite{bowman1,silverman1}.  Twenty runs for each kernel were performed using a 2.9 GHz Intel$^{\textrm{\textregistered}}$ Core i9 Macbook Pro running Mac OS-X version 10.13.6 (High Sierra) with 32 GB 2400 MHz DDR4 RAM.  The times recorded in this paper are the times for each run.    

\subsection{Soft Lunar Landing Optimal Control Problem}\label{subsect:example1}
Consider the following chance constrained modification of the soft lunar landing optimal control problem.  Minimize the cost functional
\begin{equation}\label{eq:costvers1}
\min \ J = \int_{0}^{t_f} u d t,
\end{equation} 
subject to the dynamic constraints
\begin{equation} \label{eq:dynvers1}
\begin{array}{ccc}
\dot y_1 (t) & = & y_2 (t), \\ 
\dot y_2 (t) & = & g+u(t), \\ 
\end{array}
\end{equation}
the boundary conditions
\begin{equation}\label{eq:boundvers1}
  \begin{array}{cccccc}
 \ y_1(0) & = & 10, &  y_1(t_f)& = & free, \\ 
 \ y_2(0) & = & -2, & y_2(t_f) & = & 0, \\ 
\end{array}
 \end{equation}
 the control bounds
 \begin{equation}\label{eq:contvers1}
 u \geq 0,
 \end{equation}
 the chance event constraint
 \begin{equation}\label{eq:ECC}
 \epsilon_a \geq P \left( |y_1(t_f)-\xi_1|-\delta > 0 \right),
 \end{equation}
and the chance path constraint 
\begin{equation}\label{eq:CC1}
\epsilon_b \geq P \left( u+\xi_2-3 > 0 \right),
\end{equation}
where $g = 1.622$ and $\delta=0.250$.  The risk violation parameters were set as $\epsilon_a=0.1$ and $\epsilon_b=0.01$.  The random variable $\xi_1$ was assigned a distribution of the form $N(0,0.1^2)$, while the random variable $\xi_2$ was assigned the following bimodal distribution:
\begin{equation}\label{eq:bimodal}
\xi_2 \sim \frac{1.03}{0.05 \sqrt{2 \pi}} \exp \left( \frac{-x^2}{2(0.05^2)} \right) + \frac{1.12}{0.08 \sqrt{2 \pi}} \exp \left( \frac{-(x+0.07)^2}{(0.08^2)} \right).
\end{equation}
The initial guess was a straight line approximation between the known initial and terminal conditions.  If endpoint conditions were not available, a constant initial guess that did not violate the constraint bounds was used.  The control was set as a constant of zero for the initial guess.

In order to better compare results, the solution to the following deterministic soft lunar landing optimal control problem was also computed.  Minimize Eq.~\eqref{eq:costvers1}, subject to the dynamic constraints of Eq.\eqref{eq:dynvers1}, the boundary conditions 
\begin{equation}\label{eq:boundvers2}
  \begin{array}{cccccc}
 \ y_1(0) & = & 10, &  y_1(t_f)& = & 0, \\ 
 \ y_2(0) & = & -2, & y_2(t_f) & = & 0, \\ 
\end{array}
 \end{equation} 
 and the control bounds
 \begin{equation}\label{eq:ex1path}
0 \leq u \leq 3.
\end{equation}  

 \begin{figure}[ht!]
 \centering
 \vspace*{0.25in}
\subfloat[Example problem position.]{\includegraphics[height = 2.4in]{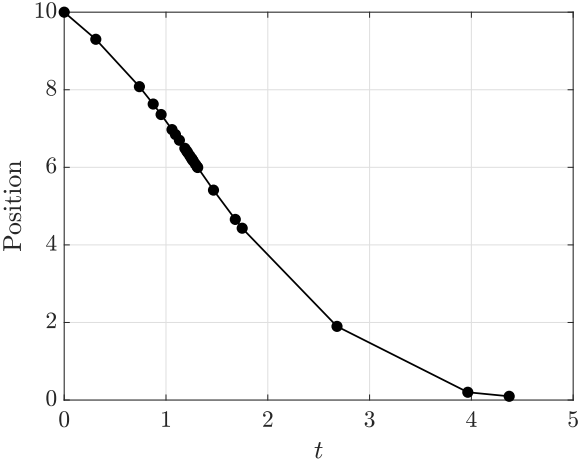}}
~~~~\subfloat[Example problem velocity.]{\includegraphics[height = 2.4in]{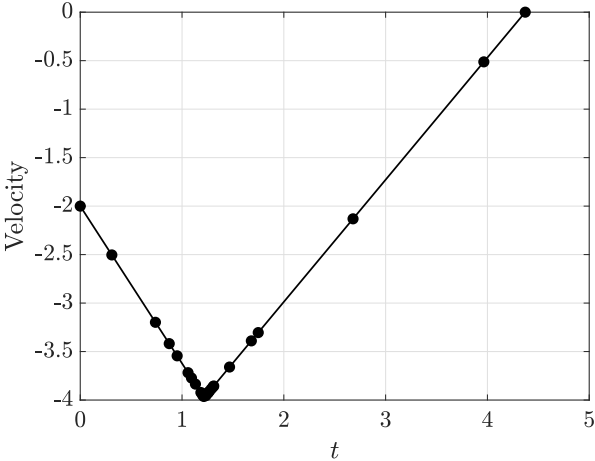}} \\
\subfloat[Example problem control.]{\includegraphics[height = 2.4in]{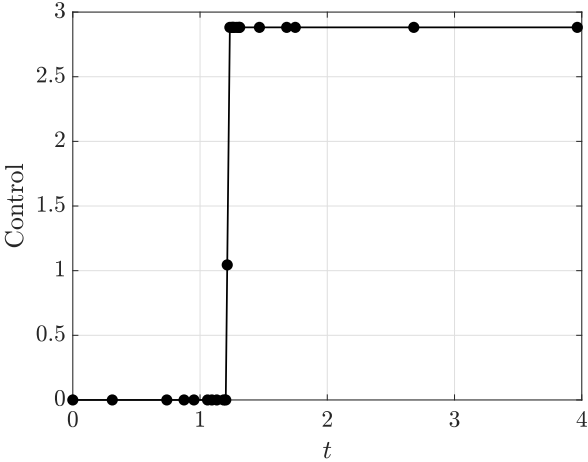}}
\caption{Example results using Split-Bernstein kernel.}
\label{fig:LunarProblem}
\end{figure}

As mentioned in Section~\ref{sect:method}, evaluation of the chance constraint at each collocation point can be computationally expensive because of the large number of MCMC samples.  In order to reduce the computation time required to solve the NLP, the chance constraint of Eq.~\eqref{eq:CC1} was replaced by the following chance constraint:  
\begin{equation}\label{eq:CCbarr}
\epsilon_b \geq 
\begin{cases}
0, \ \textrm{if} \ u \leq 3-b, \\
P \left( u+\xi_2-3 > 0 \right), \ \textrm{if} \ u > 3-b,
 \end{cases}
\end{equation}
where $b = 1$.  The reformulated chance constraint of Eq.~\eqref{eq:CCbarr} does not change the solution if $b$ is properly chosen.  In particular, looking at the form of the chance constraint in Eq.~\eqref{eq:CC1}, unless the control is close to its maximum the chance constraint will evaluate to a small number.  Consequently, implementation of the path chance constraint of Eq.~\eqref{eq:CCbarr} prevents the NLP solver from unnecessarily evaluating $50,000+$ MCMC samples at every collocation point.

The solution obtained for one run applying the Split-Bernstein kernel is shown in Fig.~\ref{fig:LunarProblem}, where it is noted that all three kernels have similar solutions and thus only one figure is included for reference.  The control still maintains the bang-bang structure of the solution to the deterministic optimal control problem with a lower maximum control than that of the solution to the deterministic optimal control, in order to satisfy the path chance constraint.  Additionally, all three kernels satisfy the event constraint as the average final positions for the Split-Bernstein, Gaussian, and Epanechnikov kernels are, respectively, $0.1100$, $0.1100$, and $0.1111$.  

\begin{table}[ht]
\caption{Run results for Example}
\renewcommand{\arraystretch}{1}
\label{table2}
\begin{center}
 	\begin{tabular}{| c || c | c | c | c |}
 	\hline
 	& Gaussian & Split-Bernstein & Epanechnikov & Deterministic \\ \hline
 	 $\mu_{J^*}$ & $9.1375$ & $9.0934$ & $9.0909$ & $8.9069$ \\ \hline
 	 $ \sigma_{J^*}$ & $0.0042$ & $0.0034$ & $0.0031$ & $0$ \\ \hline
 	 $\mu_T$ & $3.2567$ s & $4.0141$ s & $5.9531$ s & $0.1109$ s \\ \hline
 	 $\sigma_T$ & $1.1129$ s & $0.8698$ s & $1.6313$ s & $0.02039$ s \\ \hline
 	 $T_{\max}$ & $5.9796$ s & $6.3566$ s & $10.4826$ s & $0.2073$ s \\ \hline
 	 $T_{\min}$ & $1.8637$ s & $3.1013$ s & $3.6148$ s & $0.0979$ s \\ \hline 
\end{tabular}
\end{center}
\end{table} 

Table~\ref{table2} compares the results obtained using the algorithm outlined in Section~\ref{sect:method} for twenty runs of the chance constrained version of the example applying each of the three kernels, alongside the results obtained for the deterministic formulation of the example. In Table~\ref{table2}, $J^*$ is the optimal cost and $T$ is the total time for a call to the optimal control software to reach convergence for one run, where the time $T$ does not include the time to generate the initial setup for each run.  For the results shown in Table~\ref{table2}, the bandwidths for the path and event chance constraints were initially set to $0.01$ and $0.02$, respectively, for the Gaussian kernel and were set to $0.02$ and $0.03$, respectively, for the other two kernels.  Applying the criteria provided in the algorithm in Section~\ref{sect:method}, the bandwidths were switched to the bandwidths determined using $ksdensity$ for the path and event chance constraints as approximately $0.008$ and $0.01$, respectively, for all three kernels.  The results shown in Table~\ref{table2} compare well with the results obtained using the Split-Bernstein approximation method from Ref.~\cite{Kumar1}.  The average optimal costs obtained using the three kernels are larger than the optimal cost from Ref.~\cite{Kumar1}.  This discrepancy in optimal cost is due to the application of the automated bandwidth selection method resulting in bandwidths that were different for each constraint, as well as being larger than the bandwidth applied for both constraints in Ref.~\cite{Kumar1}.  The automated bandwidth selection method determined those bandwidths that most effectively balanced the trade-off between the error and variance of the biased KDE applied to each constraint.  

Next, comparing the results for the three kernels from Table~\ref{table2}, the Gaussian kernel had a lower run time than the other two kernels, as expected due to the smoothness of the function $K_G$.  The Epanechnikov kernel had the worst average run time, due to more internal operations being required when applying this kernel.  The average optimal costs between the three kernels do not vary significantly, with the average optimal cost for the Epanechnikov kernel being the smallest.  The average optimal cost for the Gaussian kernel was the largest as a result of this kernel having a larger bias than those of the other kernels.  Moreover, the run times for solving both the chance constrained and deterministic formulations of the example are small enough to be considered computationally efficient.  Finally, while the average optimal costs are higher than the deterministic optimal cost, the differences are small.    

\section{Conclusions}\label{sect:conc}
A method has been developed for approximating chance constraints as nonlinear constraints using biased KDEs with MCMC sampling to transform CCOCPs to deterministic optimal control problems.  This method combines the advantages of previous method for reformulating the chance constraints as the biased KDEs not only provide an upper bound on the chance constraint, but are also approximate CDFs.  The method of biased KDEs was then applied to an example for three different kernels.  Two of the kernels chosen, the Split-Bernstein and Epanechnikov kernels, were chosen as they satisfied the necessary requirements to obtain the correct form of a biased KDE.  The third kernel chosen was the Gaussian kernel (added for the smoothness of its integrated function), despite not meeting the qualifications to obtain the correct form of a biased KDE.  The Gaussian kernel had the lowest average computation time among the kernels, indicating that smoothness affects computational efficiency.  The Epanechnikov kernel had the lowest average optimal cost, which was still higher than the optimal cost for a deterministic formulation of the example.  Transforming CCOCPs to deterministic optimal control problems using the method of biased KDEs with MCMC sampling shows promise as the resulting optimal control problem can be solved using available optimal control software.  With fine tuning, the computational performance of this method could be made comparable to results obtained using a conservative deterministic formulation of the same CCOCP.

\section{Acknowledgments}\label{sect:ack}

The authors gratefully acknowledge support for this research from the from the U.S.~National Science Foundation under grants CMMI-1563225, DMS-1522629, and DMS-1819002.  

\renewcommand{\baselinestretch}{1.0}
\normalsize\normalfont
\bibliographystyle{aiaa}

\end{document}